\documentclass[12pt,a4paper]{amsart}
\usepackage[top=35mm, bottom=30mm, left=30mm, right=30mm]{geometry}
\usepackage{mathptmx}
\usepackage{mathrsfs}

\usepackage{verbatim}
\usepackage{url}
\usepackage[all]{xy}
\usepackage{xcolor}
\usepackage[colorlinks=true,citecolor=blue]{hyperref}

\usepackage{amsmath,amsthm}
\usepackage{amssymb}

\usepackage{enumerate}

\usepackage{graphicx}

\newtheorem{thm}{Theorem}[section]

\newtheorem{lem}[thm]{Lemma}
\newtheorem{prop}[thm]{Proposition}

\newtheorem{ques}[thm]{Question}

\theoremstyle{definition}
\newtheorem{defin}[thm]{Definition}
\newtheorem{rem}[thm]{Remark}
\newtheorem{exa}[thm]{Example}

\numberwithin{equation}{section}

\begin{document}


\baselineskip=17pt


\title{The realization and classification of topologically transitive group actions on $1$-manifolds}

\author{Enhui Shi}

\address[E.H. Shi]{School of Mathematical Sciences, Soochow University, Suzhou 215006, P. R. China}
\email{ehshi@suda.edu.cn}

\begin{abstract}
In this report, we first recall the Poincar\'e's classification theorem for minimal orientation-preserving homeomorphisms
on the circle and the Ghys' classification theorem for minimal orientation-preserving group actions on the circle. Then we
introduce a classification theorem for a specified class of topologically transitive orientation-preserving group actions on the
circle by $\mathbb Z^d$. Also, some groups that admit/admit no topologically transitive actions on the line are determined.
\end{abstract}

\maketitle

\pagestyle{myheadings} \markboth{E. H. Shi }{The realization and classification}

\section{Introduction}

The classical theory of dynamical systems study the orbit structure of actions of the group $\mathbb Z$ or $\mathbb R$. The dynamics of general group actions on
manifolds have also been extensively studied by people from different fields. One of the most basic questions in this area is:
{\it given a group $G$ and a manifold $M$, can we classify all the actions of $G$ on $M$ up to topological conjugations?} Throughout the paper, all groups
are referred  to countable discrete groups.

Certainly, the answer to this question depends on the topology of $M$, the algebraic structure of $G$, and the dynamics of the action of $G$ on $M$.
Topological transitivity and minimality are two kinds of natural irreducibility for group actions; the line and the circle are the most simple
manifolds; so, it is natural to classify  minimal/topologically transitive group actions on one dimensional manifolds first. We propose the following
question:

\begin{ques}
Given a group $G$, does it admit a minimal/topologically transitive action on 1-manifolds? If so, can we classify all such actions up to topological conjugations?
\end{ques}

The Poincar\'e's classification theorem around minimal circle homeomorphisms is the first celebrated result toward the answer to this question and the rotation numbers are complete invariants for such systems (see \cite{Po});
Ghys classified all orientation-preserving minimal group actions on the circle using bounded Euler class which extended the previous theorem due to Poincar\'e (see \cite{Gh87} or \cite{Gh01}).

Topological transitivity is weaker than minimality and the phenomena of topologically transitive group actions on the circle are much richer than that of minimal group actions.
Shi and Zhou determined some groups which admit topological transitive actions on the line (see \cite{SZ1}) and  established a classification theorem for
a specified class of topologically transitive $\mathbb Z^d$-actions on the line (see \cite{SZ2}) where $d\geq 2$. Based on the results in \cite{SZ2}, Shi and Xu further obtained a classification
theorem for tightly transitive and almost minimal (see Section 6 for the definitions) $\mathbb Z^d$-actions on the circle in \cite{SX}.
Roughly speaking, all the conjugation classes are parameterized by a combination of orbits of irrational numbers under the action of $GL(2, \mathbb Z)$ by M\"obius transformations and orbits of $\mathbb Z^n$ under some specified affine actions. This is also an extension of the Poincar\'e's classification theorem mentioned above, which indicates that all conjugation classes are parameterized by the orbits of irrationals under the $\mathbb Z$ action on $\mathbb R$ generated by the unit translation.

We compare our classification theorem with the Poincar\'e's classification theorem in the following tublar. (``TT" and ``AM" denote the properties of tight transitivity and almost minimality respectively; $\mathcal{O}(...)$ denotes the orbits of ... ).
\vspace{5mm}

\begin{tabular}{|c|c|c|}\hline
  & Poincar\'e's classification & Shi-Xu-Zhou's classification\\ \hline
Spaces & $\mathbb S^1$ & $\mathbb S^1$\\ \hline
Groups & $\mathbb Z$ & $\mathbb Z^n\ (n\geq 2)$\\ \hline
Dynamics & Minimality & TT \& AM\\ \hline
Invariants & $\mathcal{O}$(integer translations) &  $\mathcal{O}$(M\"obius \& affine actions) \\ \hline

\end{tabular}

\vspace{5mm}
The goal of this paper is to introduce the mentioned classification theorems above.  We will
give some basic definitions and examples, and the explicit statements of the main theorems.

\section{Preliminaries}

In this section, we recall some basic definitions in the theory of dynamical systems and give some  examples
to illustrate the definitions.

\begin{defin}
Let $X$ be a topological space and let ${\rm Homeo}(X)$ be
the homeomorphism group of $X$. Suppose $G$ is a group. A group
homomorphism $\phi: G\rightarrow {\rm Homeo}(X)$ is called an {\it
action} of $G$ on $X$. The action $\phi$ is said to be {\it faithful}
if it is injective. We use the pair $(X, G)$ to denote the action of $G$
on $X$.  For brevity, we usually use $gx$ or $g(x)$ instead of $\phi(g)(x)$.
\end{defin}

\begin{exa}\label{subgroup action}
If $f$ is a homeomorphism on $X$, then $n\mapsto f^n$ ($n\in\mathbb Z$)
defines an action of $\mathbb Z$ on $X$; this action is called the {\it action
generated by $f$}. If $H$ is a subgroup of ${\rm Homeo}(X)$, then the action of $H$ on $X$  always refers to
the inclusion homomorphism $i:H\rightarrow{\rm Homeo}(X)$.
\end{exa}

\begin{defin}
The {\it orbit}
of $x\in X$ under the action of $G$ is the set $Gx\equiv\{gx:g\in
G\}$; $x\in X$ is called an $n$-{\it periodic point} of $G$ if $Gx$ consists of $n$ points;
if $n=1$, then $x$ is called a {\it fixed point} of $G$.
\end{defin}

\begin{defin}
Suppose $(X, G)$ is an action of $G$ on $X$. $(X, G)$ is called {\it topologically transitive} if
$Gx$ is dense in $X$ for some $x\in X$ and $x$ is called a {\it topologically transitive point}; is called {\it minimal} if $Gx$ is dense for every $x\in X$.
\end{defin}

Clearly, minimality implies topological transitivity by the definitions.

\begin{exa}\label{minimal rotation}
Let $\mathbb S^1$ be the unit circle in the complex plane, that is $\mathbb S^1=\{e^{2\pi{\rm i} x}: x\in\mathbb R\}$.
Suppose $\alpha$ is an irrational number. Define $\rho_\alpha:\mathbb S^1\rightarrow\mathbb S^1$ by $e^{2\pi{\rm i} x}\mapsto e^{2\pi{\rm i} (x+\alpha)}$.
Then the $\mathbb Z$ action generated by $\rho_\alpha$ is minimal.
\end{exa}

\begin{exa}
Let $SL(2, \mathbb Z)$ be the group consists of matrices with integer entries and with determinant $1$. If we
view  $SL(2, \mathbb Z)$ as linear maps on the plane, then  $SL(2, \mathbb Z)$ induces a minimal action on the
projective space ${\mathbb {RP}}^1$ which is homeomorphic to the circle. Clearly, this action is not faithful, since the action
of $-I$ is identity.
\end{exa}

\begin{exa}\label{minimal example}
For $a\in \mathbb R$, define $L_a:\mathbb R\rightarrow \mathbb R$ by $L_a(x)=x+a$ for every $x\in \mathbb R$. Let $\alpha$
be irrational and let $H=\langle L_1, L_\alpha\rangle$ be the group generated by $L_1$ and $L_\alpha$. Then the action
on $\mathbb R$ by $H$ is minimal (see Example \ref{subgroup action}).
\end{exa}

\begin{exa}
We view $\mathbb S^1$ as the one point compactification of the line $\mathbb R$. Then the group $H$ in Example \ref{minimal example} induces naturally
an action on $\mathbb S^1$ with a fixed point, which is topologically transitive but not minimal.
\end{exa}

\begin{defin}
Let $\phi, \psi: G\rightarrow {\rm Homeo}(X)$ be two actions  of $G$ on $X$. The action $\phi$ is said to be {\it topologically conjugate}
to $\psi$ if there is a homeomorphism $h:X\rightarrow X$ such that $h\phi(g)=\psi(g) h$ for every $g\in G$. If $f$ and $g$ are two homeomorphisms on
$X$ and the $\mathbb Z$ action generated by $f$ is topologically conjugate to the action generated by $g$, then we call $f$ is {\it topologically conjugate}
to $g$.
\end{defin}

\begin{exa}
Let $f$ be an orientation-preserving homeomorphism on the real line $\mathbb R$ without fixed point. Then $f$ is topologically conjugate
to the unit translation $L_1$. In fact, the conjugation $h$ can be defined as follows: fix a homeomorphism $\varphi: [0, f(0))\rightarrow [0, 1)$;
define $h(x)=\varphi(f^{-n}(x))+n$ for each integer $n$ and for every $x\in [f^n(0), f^{n+1}(0))$.
\end{exa}

\section{Poincar\'e's classification}

In this section, we are going to introduce the well-known Poincar\'e's classification theorem for minimal
orientation-preserving homeomorphisms on the circle.

Let $\mathbb S^1$ be the unit circle in the complex plane and let $\pi:\mathbb R\rightarrow \mathbb S^1$ be the covering map defined by
$\pi(x)=e^{2\pi{\rm i}x}$ for every $x\in\mathbb R$.

 \begin{lem}\label{lifting existence}
 Suppose $f$ is an orientation-preserving homeomorphism on the circle $\mathbb S^1$. Then there is an orientation-preserving homeomorphism
 $F$ on the line $\mathbb R$ such that $f(\pi(x))$ $=\pi (F(x))$ for every $x\in\mathbb R$.
 \end{lem}

The map $F$ in Lemma \ref{lifting existence} is said to be a {\it lifting} of $f$.

 \begin{lem}\label{lifting lemma}
Let $F$ and $G$ be two liftings of $f$. Then there is a unique integer $n$ such that $F(x)=G(x)+n$ for every $x\in\mathbb R$.
 \end{lem}

 \begin{prop}
Let $f$ be an orientation-preserving homeomorphism on $\mathbb S^1$ and let $F:\mathbb R\rightarrow \mathbb R$ be a lifting of $f$.
Fix an $x\in\mathbb R$. Then $\lim\limits_{n\to\infty}\frac{F^n(x)-x}{n}$ exists and is independent of the choice of $x$.
 \end{prop}

If we denote this limit by $\alpha$, then $\alpha\ {\rm mod}\ 1$ is independent of the choice of $F$ from Lemma \ref{lifting lemma}.

\begin{defin}
The number $\alpha\ {\rm mod}\ 1$  is said to be the {\it rotation number} of $f$, which is denoted by $\rho(f)$.
\end{defin}

\begin{exa}
The rotation number of $\rho_\alpha$ in Example \ref{minimal rotation} is $\alpha\ {\rm mod}\ 1$.
\end{exa}

The following proposition shows that rotation numbers are conjugation invariants for orientation-preserving homeomorphisms on
the circle.

\begin{prop}
Let $f$ and $g$ be two orientation-preserving homeomorphisms on the circle. If $f$ is topologically conjugate to $g$, then $\rho(f)=\rho(g)$.
\end{prop}

The following classification theorem is due to H. Poincar\'e (see \cite{Po}).

\begin{thm}
Let $f$ be an orientation-preserving homeomorphisms on the circle which is minimal. Then the rotation number $\rho(f)$ is
irrational and $f$ is topologically conjugate to $\rho_\alpha$ as in Example \ref{minimal rotation}, where $\alpha=\rho(f)$.
\end{thm}

\section{Ghys' classification}

For orientation-preserving group actions on $\mathbb S^1$, Ghys established a classification theorem in \cite{Gh87}, which
is an analogy of the Poicar\'e's classification theorem. The invariant used in Ghys' classification theorem is bounded Euler class.
Now we will recall some notions before the statement of Ghys' Theorem. All the contents in this section come from \cite{Gh01}.

\begin{defin}
Let $\Gamma$ be a group and let $A$ be an abelian group. A {\it k-cochain} of $\Gamma$ with values in $A$ is a map $c:\Gamma^{k+1}\rightarrow A$
that is {\it homogeneous} (that is, $c(\gamma\gamma_0, \gamma\gamma_1, ..., \gamma\gamma_k)$=$c(\gamma_0, \gamma_1, ...,\gamma_k)$ indentically).
The set of $k$-cochains is denoted by $C^k(\Gamma, A)$, which is an abelian group. The {\it coboundary map} $d_k:C^k(\Gamma, A)\rightarrow C^{k+1}(\Gamma, A)$
is defined by
$$
d_kc(\gamma_0,...,\gamma_{k+1})=\sum_{i=0}^{k+1}(-1)^ic(\gamma_0,...,\hat\gamma_i,...,\gamma_{k+1}).
$$
Then $d_{k+1}\circ d_k=0$ and the {\it cohomology group} $H^k(\Gamma, A)$ is defined to be the quotient ${\rm Ker}(d_k)/{\rm Im}(d_{k-1})$. Elements in
${\rm Ker}(d_k)$ are called {\it cocycles} and elements in ${\rm Im}(d_{k-1})$ are called {\it coboundaries}.
\end{defin}

\begin{lem}
A map $c:\Gamma^{k+1}\rightarrow A$ is homogeneous if and only if there is a (unique) $\bar c:\Gamma^k\rightarrow A$ such that
$c(\gamma_0,...,\gamma_k)={\bar c}(\gamma_0^{-1}\gamma_1,\gamma_1^{-1}\gamma_2,...,\gamma_{k-1}^{-1}\gamma_k)$.
\end{lem}

\begin{defin}
The map $\bar c$ in the above lemma is called the {\it inhomogeneous cochain} associated to $c$.
\end{defin}

We consider a central extension of $\Gamma$ by $A:0\rightarrow A \stackrel{i}{\rightarrow} \tilde \Gamma \stackrel{p}\rightarrow \Gamma \rightarrow 1$,
which means that the center of $\tilde \Gamma$ contains a subgroup isomorphic to $A$ and the quotient of $\tilde \Gamma$ by this subgroup is isomorphic
to $\Gamma$. Choose a set theoretical section $s$ from $\Gamma$ to $\tilde \Gamma$ and let $\bar c(\gamma_1, \gamma_2)=s(\gamma_1\gamma_2)^{-1}s(\gamma_1)s(\gamma_2)$.
Then $p(\bar c(\gamma_1, \gamma_2))$ is the identity in $\Gamma$ and can be identified with an element of $A$. Let $c:\Gamma^3\rightarrow A$ be the associated homogeneous cochain of
$\bar c:\Gamma^2\rightarrow A$ (see Lemma 4.2). One can check that $c$ is a cocycle.

\begin{defin}
The cohomology class of $c$ in $H^2(\Gamma, A)$ is called the {\it Euler class of the extension} under consideration (which does not depend on the choice of a section).
\end{defin}

Let $\widetilde {\rm Homeo}_+(\mathbb S^1)$ be the group of homeomorphisms of $\mathbb R$ which commute with integral transformations and let ${\rm Homeo}_+(\mathbb S^1)$
be the group of all orientation-preserving homeomorphisms of $\mathbb S^1$.
Then $\widetilde {\rm Homeo}_+(\mathbb S^1)$
consists of all liftings of elements in ${\rm Homeo}_+(\mathbb S^1)$ by Lemma 3.1. Then we get the following center extension:
\begin{equation}
0\rightarrow \mathbb Z\rightarrow\widetilde {\rm Homeo}_+(\mathbb S^1)\stackrel{p}{\rightarrow}{\rm Homeo}_+(\mathbb S^1)\rightarrow 1.
\end{equation}

Let $\phi:\Gamma\rightarrow {\rm Homeo}_+(\mathbb S^1)$ be an action of $\Gamma$ on $\mathbb S^1$ and let
$\tilde \Gamma=\{(\gamma, \tilde f)\in\Gamma\times\widetilde {\rm Homeo}_+(\mathbb S^1):\phi(\gamma)=p(\tilde f)\}$.
Let $p':\tilde \Gamma\rightarrow \Gamma$ be the canonical projection. Then we have a central extension of $\Gamma$ by $\mathbb Z$:
\begin{equation}
0\rightarrow \mathbb Z\rightarrow\tilde \Gamma \stackrel{p'}{\rightarrow}\Gamma\rightarrow 1.
\end{equation}

\begin{defin}
The Euler class corresponding to the extension $(4.2)$ is called the {\it Euler class of the action $\phi$} and is denoted by $\phi^*(eu)$.
\end{defin}

Though $\phi^*(eu)$ is a conjugation invariant, it is very poor, since it cannot even defect the rotation number by the fact that $H^2(\mathbb Z, \mathbb Z)=0$.
 Ghys used bounded Euler classes to fill up this deficiency.

\begin{defin}
Let $A=\mathbb Z$ or $\mathbb R$. A bounded homogeneous map from $\Gamma^{k+1}$ to $A$ is called a {\it bounded $k$-cochain}.
\end{defin}

The set of all bounded $k$-cochain is denoted by $C_b^k(\Gamma, A)$, which is a sub $A$-module of $C^k(\Gamma, A)$.
It is clear that the coboundary of $d_k$ of a bounded $k$-cochain is a bounded $(k+1)$-cochain.

\begin{defin}
The cohomology of this new differential complex is called the {\it bounded cohomology of $\Gamma$ with coefficients in A} and is denoted by $H_b^k(\Gamma, A)$.
\end{defin}

Now, let us look at $(4.1)$. For every $f\in {\rm Homeo}_+(\mathbb S^1)$, let the unique $\sigma(f)\in p^{-1}(f)$ be such that $\sigma(f)(0)\in[0,1)$.
Then $\sigma:{\rm Homeo}_+(\mathbb S^1)\rightarrow \widetilde {\rm Homeo}_+(\mathbb S^1)$ is a section. Let $\bar c$ be the associated inhomogeneous cocycle,
that is $\bar c(f_1,f_2)=\sigma(f_1f_2)^{-1}\sigma(f_1)\sigma(f_2)$. It is easy to see that the associated $2$-cocycle $c$ is bounded and integral.

\begin{defin}
The above element $c\in H_b^2({\rm Homeo}_+(\mathbb S^1), \mathbb Z)$ is called the {\it bounded Euler class}. If $\phi:\Gamma\rightarrow {\rm Homeo}_+(\mathbb S^1)$
is an action, then the element $\phi^*(eu)$ in $H_b^2(\Gamma, \mathbb Z)$  is called the {\it bounded Euler class} of $\phi$.
\end{defin}

Consider a group homomorphism $r:\Gamma\rightarrow \mathbb Z/k\mathbb Z$ where $k$ is a positive integer. Then $r$ naturally induces an action of
$\Gamma$ on $\mathbb S^1$ by rotations of order $k$.

\begin{defin}
The bounded Euler class of the action induced by $r$ is called a {\it rational element} in $H_b^2(\Gamma, \mathbb Z)$.
\end{defin}

The following classification theorem is due to Ghys. We only take a part of the main theorems in \cite{Gh87}, since we are only interested in
topologically transitive actions in this report.

\begin{thm}
Suppose $\phi_1,\phi_2:\Gamma\rightarrow {\rm Homeo}_+(\mathbb S^1)$ are two minimal actions of $\Gamma$ on $\mathbb S^1$. Then $\phi_1$
is topologically conjugate to $\phi_2$ if and only if $\phi_1^*(eu)=\phi_2^*(eu)$ (which is not rational).
\end{thm}

\section{The realization of topologically transitive group actions}

In this section, we are mainly concerned with the existence of topologically transitive group actions on the line $\mathbb R$. There are two ways that group actions on $\mathbb R$
are related to group actions on $\mathbb S^1$. First, every orientation-preserving action of a group $G$ on $\mathbb R$ naturally induces an orientation-preserving action on the circle $\mathbb S^1$
with a fixed point by two points compactification of $\mathbb R$. Second, every orientation-preserving action of a group $G$ on $\mathbb S^1$ can be lifted to a group action on $\mathbb R$ by central extension
(see $(4.1)$). So, we should first study group actions on the line.

\begin{defin}
Let $\phi:G\rightarrow {\rm Homeo}_+(\mathbb R)$ be an action of group $G$ on the line.  An open interval $(a,b)\subset \mathbb R$ is said to be
a {\it wandering interval} of $\phi$ or of $G$ if, for every $g\in G$, either the restriction $g|_{(a,b)}={\rm Id}_{(a,b)}$
or $g((a, b))\cap (a, b)=\emptyset$.
\end{defin}

\begin{exa}
Every interval $(a, b)$ with $|b-a|<1$ is wandering for the $\mathbb Z$ action on $\mathbb R$ generated by the unit translation $L_1$.
\end{exa}

Shi and Zhou obtained the following dichotomy theorem, which is useful in determining whether a group admits a topologically
transitive action on $\mathbb R$ (see \cite{SZ1}).

\begin{thm}\label{dichonomy thm}
Let $G$ be a group. Then either $G$ has a topologically transitive action on the line $\mathbb R$ by orientation-preserving homeomorphisms,
or every orientation-preserving action of $G$ on $\mathbb R$ has a wandering interval.
\end{thm}

We should note that the ``dichotomy phenomenon" in Theorem \ref{dichonomy thm} is far from being true for group actions on spaces of dimension $\geq 2$.
For example, if $D$ is the closed unit disk in the plane and $S^2$ is the unit sphere in $\mathbb R^3$, then every one-point union of $D$ and $S^2$
admits no topologically transitive homeomorphism but admits a homeomorphism with no wandering open set.

\begin{defin}
A group $G$ is {\it poly-cyclic} (resp. {\it super-poly-cyclic}) if it admits a decreasing sequence of subgroups
$G=N_0\rhd N_1\rhd...\rhd N_k=\{e\}$ for some positive integer $k$ such that $N_{i+1}$ is normal in $N_i$ (resp. $N_{i+1}$ is normal in $G$) and $N_i/N_{i+1}$ is cyclic for each $i\leq k-1$; it is called {\it poly-infinite-cyclic} (resp. {\it super-poly-infinite-cyclic}) if $N_i/N_{i+1}$ is infinitely cyclic for each $i\leq k-1$.
\end{defin}

It is well known that all poly-cyclic groups are solvable and all finitely generated torsion free nilpotent groups are super-poly-infinite-cyclic.

\begin{defin}
Suppose $G=N_0\rhd N_1\rhd...\rhd N_k=\{e\}$ is super-poly-infinite-cyclic. Take $f_i\in N_i\setminus N_{i+1}$ such that $N_i/N_{i+1}=\langle f_iN_{i+1}\rangle$ for each $i$. Then $f_if_{i+1}N_{i+2}f_i^{-1}=f_{i+1}^{n_i}N_{i+2}$ where $n_i=\pm 1$. We call the $(k-1)$-tuple $(n_0, n_1, ..., n_{k-2})$ {\it the name of $G$}. Clearly, the name of $G$ is independent of the choice of $f_i$.
\end{defin}

\begin{defin}
For each integer $n$, the solvable {\it Baumslag-Solitar group} $B(1, n)$ is the group $\langle a, b:  ba=a^nb \rangle$.
\end{defin}

The Baumslag-Solitar groups are examples of two-generator one-relator groups that play an important role in combinatorial group theory and geometric group theory.
Notice that $B(1, 1)$ is the free abelian group of rank $2$ and $B(1, -1)$ is the fundamental group of the  Klein Bottle, which is a classical example
being of orderable but non bi-orderable (see \cite[Exercise 2.2.68]{Na1}).

Based on Theorem \ref{dichonomy thm}, Shi and Zhou get the following theorem in \cite{SZ1}.

\begin{thm}\label{existence thm}
The following groups admit topologically transitive orientation-preserving actions on $\mathbb R$: the nonabelian free group $\mathbb Z*\mathbb Z$; any super-poly-infinite-cyclic group $G=N_0\rhd N_1\rhd...\rhd N_k=\{e\}$ having the name $(n_0, n_1, ..., n_{k-2})$ with some $n_i=1$ and with $k\geq 2$; any poly-infinite-cyclic, non super-poly-infinite-cyclic group $G$; the Baumslag-Solitar group $B(1, n)$ with $n\not=0$ and $n\not=-1$.

The following groups admit no topologically transitive orientation-preserving actions on $\mathbb R$: finite groups;  the infinite cyclic group $\mathbb Z$; $SL(2, \mathbb Z)$; finite index subgroups of $SL(n,\mathbb Z)$ with $n\geq 3$; any super-poly-infinite-cyclic group $G=N_0\rhd N_1\rhd...\rhd N_k=\{e\}$ having the name $(-1, -1, ..., -1)$; the Baumslag-Solitar group $B(1, -1)$.
\end{thm}

\begin{rem}
The nonexistence of topologically transitive orientation-preserving actions on $\mathbb R$ for finite index subgroups of $SL(n,\mathbb Z)$ with $n\geq 3$ in the above theorem is almost a restatement of the celebrated  result
due to Witte-Morris in \cite{Wi2}.
\end{rem}

From Theorem \ref{existence thm}, we see that many solvable groups admit topologically transitive actions on $\mathbb R$ (and then on $\mathbb S^1$). Recall that
a {\it higher rank lattice} is a lattice of a simple Lie group with finite center and with real rank $\geq 2$.
Contrary to solvable group actions, it is conjectured that  every continuous action of a higher rank lattice on the
circle $\mathbb S^1$ must factor through a finite group action. Though the conjecture is still open now,
Ghys (see \cite{Gh99}) and Burger-Monod (see \cite{Bu}) proved independently the existence of periodic points for such actions.
This converts the conjecture to that  no higher rank lattice admits faithful actions on the line $\mathbb R$.
We get the following proposition in \cite{SZ1}.
\begin{prop}
Suppose $G$ is a higher rank lattice. If $G$ admits an orientation-preserving faithful action on $\mathbb R$,
then it also admits a topologically transitive action on $\mathbb R$.
\end{prop}

\section{The classification of topologically transitive $\mathbb Z^n$ actions}

In this section, we introduce the classification theorem for a specified class of topologically transitive $\mathbb Z^n$ actions on $\mathbb S^1$ obtained by
Shi, Xu, and Zhou in \cite{SZ2, SX}. This is an extension of Poincar\'e's classification  for minimal homeomorphisms on the circle.

\begin{defin}
Let $X=\mathbb R$ or $\mathbb S^1$. A subgroup $G$ of  $\text{Homeo}_{+}(X)$  is  {\it tightly transitive} if its action on $X$ is topologically transitive and no subgroup $H$ of $G$ with $[G: H]=\infty$ has this property; is {\it almost minimal} if it has at most countably many nontransitive points.
\end{defin}

It is clear that any cyclic group generated by an irrational rotation on $\mathbb S^1$ is tightly transitive and almost minimal.
Now we give a class of tightly transitive and almost minimal subgroups of $\text{Homeo}_{+}(\mathbb R)$, which are isomorphic to $\mathbb Z^n$ with $n\geq 2$.

\begin{exa}\label{tightly transitive on R}
 Let $\alpha$ be an irrational number in $(0,1)$ and $n\geq 2$ be an integer. Let $a,b\in\mathbb{R}$.  Denote by $\langle L_a, L_b\rangle$ the subgroup of $\text{Homeo}_{+}(\mathbb{R})$ generated by  $L_a$ and $L_b$ .

We define $G_{\alpha,n}$ inductively. Let $G_{\alpha,2}=\langle L_1, L_\alpha\rangle$. Suppose that we have constructed $G_{\alpha,n}$ for $n\geq 2$. Then
we construct $G_{\alpha, n+1}$ as follows. Let $\hbar$ be the homeomorphism from $\mathbb{R}$ to $(0,1)$ defined by
\begin{displaymath}
\hbar(x)=\frac{1}{\pi}\left(\arctan x+\frac{\pi}{2}\right)~~~\text{for } ~x\in\mathbb{R}.
\end{displaymath}
For $\sigma\in G_{\alpha,n}$, define $\hat{\sigma}\in\text{Homeo}_{+}(\mathbb{R})$ by
\begin{equation}\label{definition of hat}
\hat{\sigma}(x)=\left\{
\begin{array}{cl}
\hbar \sigma \hbar^{-1}(x-i)+i,& x\in (i,i+1)~~\text{and} ~i\in\mathbb{Z},\\
x,& x\in\mathbb{Z}.
\end{array}
\right.
\end{equation}
\noindent Let $G_{\alpha,n+1}$ be the group generated by $\{\hat{\sigma}: \sigma\in G_{\alpha,n}\}\cup\{L_1\}$. Then $G_{\alpha,n}$ is isomorphic to $\mathbb Z^n$ and is
tightly transitive and almost minimal; any isomorphism $\phi:\mathbb Z^n\rightarrow G_{\alpha,n}$ gives the desired action.

\begin{figure}[htp]
\centering
\includegraphics[width=13cm]{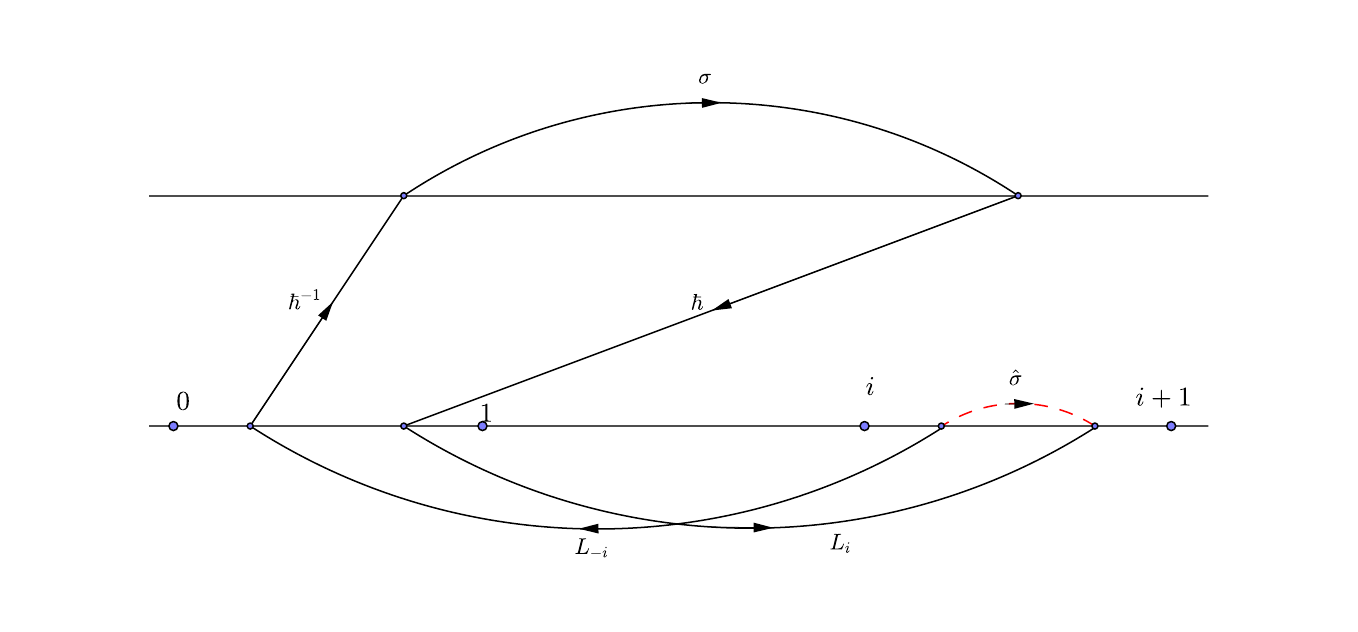}\\
\caption{Definition of $\hat{\sigma}$}
\end{figure}

\end{exa}

Suppose that $\alpha$ and $\beta$ are irrationals in $(0,1)$. We say that $\alpha$ is equivalent to $\beta$ if there exist $m_1,n_1,m_2,n_2\in \mathbb{Z}$ with $\left|m_1n_2-n_1m_2\right|=1$ such that $\beta=\frac{m_1+n_1\alpha}{m_2+n_2\alpha}$. The following classification theorem is due to Shi and Zhou (see \cite{SZ2}).

\begin{thm}\label{classification of R} The following assertions hold:

\begin{itemize}
\item[(1)] For any $n\geq 2$ and  irrationals $\alpha,\beta\in(0,1)$,  $G_{\alpha,n}$ is  conjugate to $G_{\beta,n}$ as subgroups of
 $\text{Homeo}_{+}(\mathbb{R})$ if and only if $\alpha$ is equivalent to $\beta$.
\item[(2)] Let $G$ be a subgroup of $\text{Homeo}_{+}(\mathbb{R})$ which is isomorphic to $\mathbb{Z}^n$ for some $n\geq 2$ and is tightly transitive and almost minimal. Then  $G$  is conjugate to  $G_{\alpha,n}$ for some irrational $\alpha\in(0,1)$ as subgroups of
 $\text{Homeo}_{+}(\mathbb{R})$.
\end{itemize}

\end{thm}

\begin{rem}
Let $X=\mathbb R$ or $\mathbb S^1$, and let $\phi, \psi:\mathbb Z^n\rightarrow \text{Homeo}_{+}(\mathbb{X})$ be two actions. If  $\phi$ and $\psi$ are topologically conjugate, then $\phi(\mathbb Z^n)$ and $\psi(\mathbb Z^n)$ are conjugate subgroups of $\text{Homeo}_{+}(\mathbb{X})$; conversely, if
 $\phi(\mathbb Z^n)$ and $\psi(\mathbb Z^n)$ are conjugate subgroups of $\text{Homeo}_{+}(\mathbb{X})$, then there is an automorphism $\tau:\mathbb Z^n\rightarrow \mathbb Z^n$ such that the actions $\phi$ and $\psi\circ\tau$ are topologically conjugate. So, classifying faithful $\mathbb Z^n$ actions on $X$ up to automorphisms of $\mathbb Z^n$ is equivalent to classifying the conjugation subgroups of $\text{Homeo}_{+}(\mathbb{X})$ isomorphic to $\mathbb Z^n$.
\end{rem}

Based on the examples in Example \ref{tightly transitive on R}, we now construct a class of tightly transitive and almost minimal subgroups of  $\text{Homeo}_{+}(\mathbb S^1)$, which are isomorphic to $\mathbb Z^n$ with $n\geq 2$.

\begin{exa}
Let integers  $n\geq 2$ and $k\geq 1$ and let $\alpha$ be an irrational number in $(0,1)$. Let $G_{\alpha,n}^k=\{g^k: g\in G_{\alpha,n}\}$. Suppose  $g\in G_{\alpha, n}\setminus G_{\alpha,n}^k$.  Put $x_j=e^{{\rm\bf i}2\pi j/k}$ for $j=1,...,k$. Denote by $(x_i, x_{i+1})~(\text{resp. }[x_i,x_{i+1}])$ the open (resp. closed) interval in $\mathbb S^1$ from $x_i$ to $x_{i+1}$ anticlockwise. Fix an orientation-preserving homeomorphism $\phi$ from $\mathbb{R}$ to $(x_1, x_2)$. For any $\sigma\in G_{\alpha,n}$, define $\tilde{\sigma}\in\text{Homeo}_{+}((x_1,x_2))$ by
\begin{displaymath}
\tilde{\sigma}=\phi \sigma\phi^{-1}.
\end{displaymath}
Let $f\in \text{Homeo}_{+}(\mathbb{S}^1)$ be such that
\begin{itemize}
  \item $f(x_i)=x_{i+1 }$,  for $i=1,...,k$.
  \item $f^{k}\mid_{(x_{1},x_2)}=\tilde{g}: (x_{1},x_2)\rightarrow (x_1,x_2)$.
\end{itemize}
In the above definition, we take $x_{k+1}=x_1$.  We denote the collection of such $f\in \text{Homeo}_{+}(\mathbb{S}^1)$ by $\text{Homeo}_{+}(\mathbb{S}^1)_{k,g}$.\\

Now for $\sigma\in G_{\alpha, n}$, extend $\tilde{\sigma}$ to an orientation-preserving homeomorphism $\overline{\sigma_f}$ of $\mathbb{S}^1$ by
\begin{equation}\label{extension}
\overline{\sigma_f}(x)=\left\{
\begin{array}{cll}
f^{i-1}\tilde{\sigma}f^{-(i-1)}(x),& x\in (x_i,x_{i+1}),& i=1,..., k,\\
x_{i},& x=x_i, & i=1,...,k.
\end{array}
\right.
\end{equation}

 Now we define $G_{\alpha,n,k,g,f}$ to be the subgroup of $\text{Homeo}_{+}(\mathbb{S}^1)$ generated by $\{\overline{\sigma_f}:\sigma\in G_{\alpha, n}\}\cup\{f\} $.\\
\begin{figure}[htp]
\centering
\includegraphics[width=8cm]{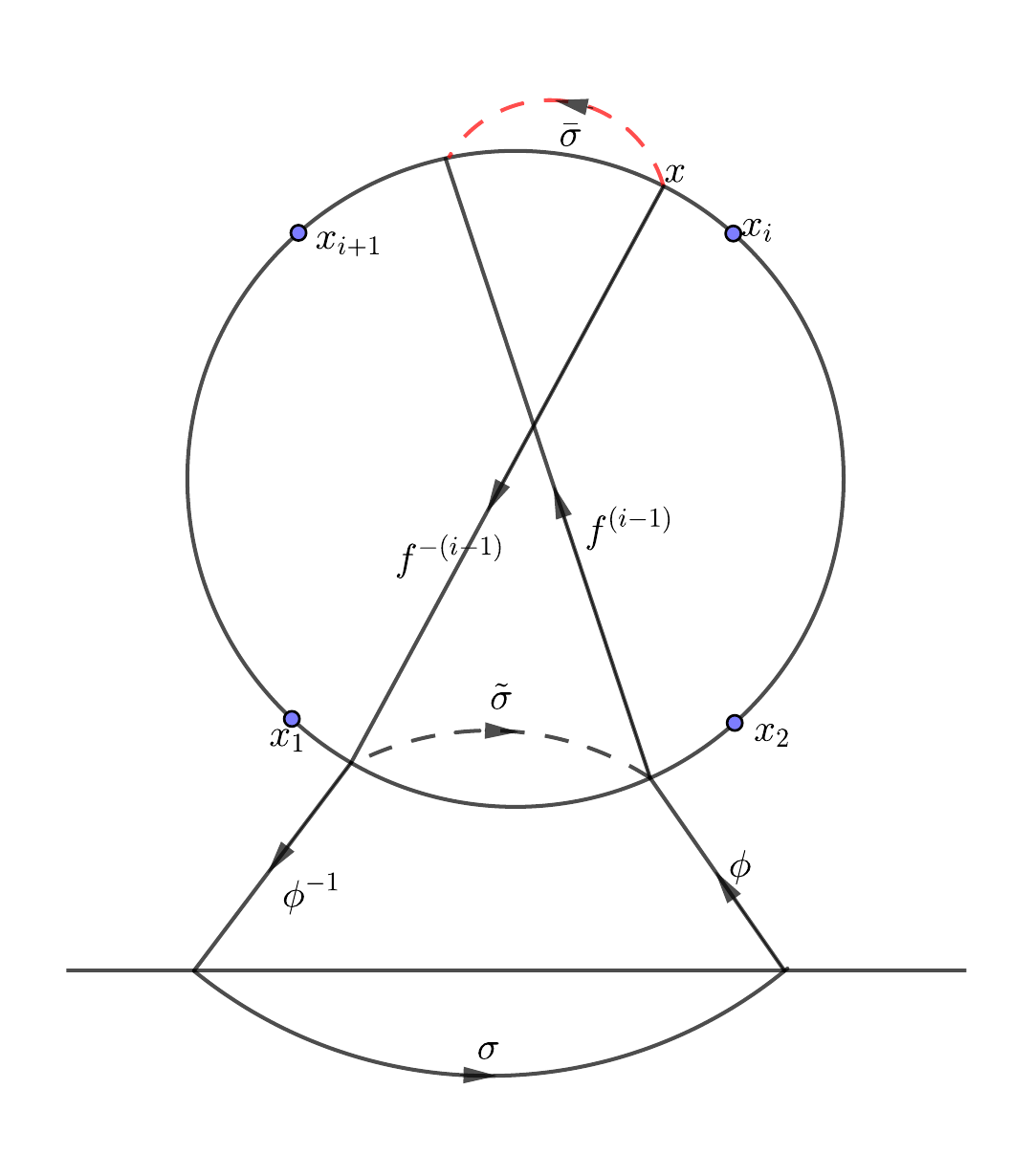}\\
\caption{Definition of $\overline{\sigma_f}$}
\end{figure}
\end{exa}
Then $G_{\alpha,n,k,g,f}$ is isomorphic to $\mathbb{Z}^n$ and  is tightly transitive and almost minimal.

The following two theorems due to Shi and Xu completely classified the topologically transitive and almost minimal
subgroups of $\text{Homeo}_{+}(\mathbb{S}^1)$, which are isomorphic to $\mathbb Z^n$ with $n\geq 2$  (see \cite{SX}).

\begin{thm}\label{all systems}
Let $G$ be a  subgroup of $\text{Homeo}_{+}(\mathbb{S}^1) $, which is isomorphic to $\mathbb{Z}^n$ for some $n\geq 2$ and is tightly transitive and almost minimal. Then $G$  is  conjugate to some $G_{\alpha,n,k,g,f}$.
\end{thm}

Let $N_{\text{Homeo}_{+}(\mathbb{R})}(G_{\alpha,n })$ denote the normalizer of $G_{\alpha,n}$ in $\text{Homeo}_{+}(\mathbb{R}) $, i.e.,
$N_{\text{Homeo}_{+}(\mathbb{R})}(G_{\alpha,n })=\left\{\varphi\in \text{Homeo}_{+}(\mathbb{R}):~\varphi G_{\alpha,n}\varphi^{-1}=G_{\alpha,n}\right\}.$
 Thus we get an affine action on $G_{\alpha,n}$ by the semidirect $N_{\text{Homeo}_{+}(\mathbb{R})}(G_{\alpha,n })\ltimes G_{\alpha,n}^k$: $(\varphi, f). g:=\varphi g\varphi^{-1}f$, for any $(\varphi, f)\in N_{\text{Homeo}_{+}(\mathbb{R})}(G_{\alpha,n })\ltimes G_{\alpha,n}^k$ and $g\in G_{\alpha,n}$.

Define
$${\rm Conj}(G_{\alpha, n}, G_{\alpha', n})=\{\psi\in{\rm Homeo}_{+}(\mathbb R):G_{\alpha, n}=\psi G_{\alpha', n}\psi^{-1}\}.$$
If $\alpha$ and $\alpha'$ are equivalent, then ${\rm Conj}(G_{\alpha, n}, G_{\alpha', n})\not=\emptyset$ by Theorem \ref{classification of R};
and we fix a conjugation $\psi_{\alpha,\alpha'}\in{\rm Conj}(G_{\alpha, n}, G_{\alpha', n})$.

\begin{thm}\label{classification of S}
 The group $G_{\alpha,n,k,g,f}$ is conjugate to $G_{\alpha',n',k',g',f'}$ if and only if
 \begin{itemize}
   \item $n=n'$ and $k=k'$;
   \item $\alpha$ is equivalent to $\alpha'$, i.e. there exist $m_1,n_1,m_2,n_2\in \mathbb{Z}$ with $\left|m_1n_2-n_1m_2\right|=1$ such that $\alpha'=\frac{m_1+n_1\alpha}{m_2+n_2\alpha}$;
   \item $g$ and $\psi_{\alpha,\alpha'} g'\psi_{\alpha,\alpha'}^{-1}$ are in the same orbit of the affine action on $G_{\alpha,n}$ by $N_{\text{Homeo}_{+}(\mathbb{R})}$ $(G_{\alpha,n })\ltimes G_{\alpha,n}^k$, i.e. there exist some $\varphi\in N_{\text{Homeo}_{+}(\mathbb{R})}(G_{\alpha,n })$ and $h\in G_{\alpha,n}$ such that $\psi_{\alpha,\alpha'} g'\psi_{\alpha,\alpha'}^{-1}=\varphi gh^k\varphi^{-1}$.
 \end{itemize}
\end{thm}

\subsection*{Acknowledgements}
The work is supported by NSFC (No. 11771318, 11790274).


\end{document}